# Optimal Distributed Energy Resources Sizing for Commercial Building Hybrid Microgrids

Yishen Wang*, Zhehan Yi*, Di Shi*, Zhe Yu*, Bibin Huang†, Zhiwei Wang*
* GEIRI North America, San Jose, CA, USA
† State Grid Energy Research Institute, Beijing, China.
Email: yishen.wang@geirina.net

*Abstract*—As microgrids have advanced from early prototypes to relatively mature technologies, converting data center integrated commercial buildings to microgrids provides economic, reliability and resiliency enhancements for the building owners. Thus, microgrid design and economically sizing distributed energy resources (DER) are becoming more demanding to gain widespread microgrids commercial viability. In this paper, an optimal DER sizing formulation for a hybrid AC/DC microgrid configuration has been proposed to leverage all benefits that AC or DC microgrid could solely contribute. Energy storage (ES), photovoltaics (PV) and power electronics devices are coordinately sized for economic grid-connected and reliable islanded operations. Time-of-use (TOU) energy usages charges and peak demand charges are explicitly modeled to achieve maximum level of cost savings. Numerical results obtained from a real commercial building load demonstrate the benefits of the proposed approach and the importance of jointly sizing DER for the grid-connected and islanded modes.

## Nomenclature

### A. Sets and Indices

$S$      Set of representative days, indexed by $s$.
$T$      Set of time intervals, indexed by $t$.

### B. Decision variables

$c^{inv}$      Total investments.
$c_s^{e/d}$      Energy/demand charges at day $s$.
$c_s^{lcl/lnl}$      Lost of critical/non-critical load costs at day $s$.
$c_s^{deg}$      Energy storage degradation costs at day $s$.
$x^{PV}$      Installed PV capacities.
$x^{ES}$      Installed batteries power ratings.
$x^{IC}$      Installed interfacing converter capacities.
$x^{INV}$      Installed inverter capacities.
$x^{CON}$      Installed converter capacities.
$p_{s,t}^{grid}$      Energy purchased from utility at day $s$ time $t$.
$p_s^{peak}$      Net load peak demand at day $s$.
$lcl_{s,t}^{AC/DC}$      Critical AC/DC load shedding at day $s$ time $t$.
$lnl_{s,t}^{AC/DC}$      Non-critical AC/DC load shedding at day $s$ time $t$.
$v_{s,t}^{DC}$      PV DC output at day $s$ time $t$.
$dch_{s,t}^{AC/DC}$      ES AC/DC discharge at day $s$ time $t$.
$ch_{s,t}^{AC/DC}$      ES AC/DC charge at day $s$ time $t$.
$soc_{s,t}$      ES state-of-charge at day $s$ time $t$.
$f_{s,t}^{AC}$      AC bus flow at day $s$ time $t$.
$f_{s,t}^{DCin}$      DC bus injection at day $s$ time $t$.
$f_{s,t}^{DCout}$      DC bus extraction at day $s$ time $t$.
$z_{s,t}^{DC}$      Binary indicator for DC bus flow direction at day $s$ time $t$, 1 for injection, 0 otherwise.
$y_{s,t}^{ES}$      Binary indicator for ES discharging status at day $s$ time $t$, 1 for discharging, 0 otherwise.
$vi_{s,t}^{DC}$      Islanded PV DC output at day $s$ time $t$.
$dchi_{s,t}^{AC/DC}$      Islanded ES AC/DC discharge at day $s$ time $t$.
$fi_{s,t}^{AC}$      Islanded AC bus flow at day $s$ time $t$.
$fi_{s,t}^{DCin}$      Islanded DC bus injection at day $s$ time $t$.
$fi_{s,t}^{DCout}$      Islanded DC bus extraction at day $s$ time $t$.
$zi_{s,t}^{DC}$      Islanded DC bus flow direction indicator at day $s$ time $t$, 1 for injection, 0 otherwise.
$u_{s,t}^{ES}/\kappa_{s,t}^{ES}$      Auxiliary variables for ES discharging state.

### C. Parameters

$\pi_s$      Probability of representative day $s$.
$C^{PV}$      Annualized investment for PV.
$C^{ES}$      Annualized investment for ES.
$C^{IC}$      Annualized investment for interfacing converters.
$C^{INV}$      Annualized investment for AC/DC inverters.
$C^{CON}$      Annualized investment for DC/DC converters.
$C^{deg}$      Degradation costs for ES.
$\lambda_t^e$      Price for energy usage charge at time $t$.
$\lambda^d$      Price for demand charge.
$VOLL^{CL}$      Value-of-lost-load for critical load.
$VOLL^{NL}$      Value-of-lost-load for non-critical load.
$P^{Tariff}$      Tariff allowed maximum peak demand.
$\eta^{IC}$      Efficiency for interfacing converter.
$\eta^{INV}$      Efficiency for DC/AC inverter.
$\eta^{CON}$      Efficiency for DC/DC converter.
$\eta^{ch/dch}$      Efficiency for ES charging/discharging.
$CL_{s,t}^{AC/DC}$      Critical AC/DC load at day $s$ time $t$.
$NL_{s,t}^{AC/DC}$      Non-critical AC/DC load at day $s$ time $t$.
$ES^{max}$      Maximum ES power capacities.
$PV^{max}$      Maximum available PV capacities.
$V_{s,t}$      PV availability profiles in p.u. at day $s$ time $t$.
$\rho^{EP}$      Batteries energy-power ratio.
$\alpha^{min/max}$      Minimum/maximum battery energy level.
$M$      Large enough constant.

## I. Introduction

A microgrid can be defined as a single controllable entity which includes a group of interconnected loads and distributed energy resources (DER) to act with respect to the grid [1].

This work is funded by SGCC Science and Technology Program.

Additionally, it can connect and disconnect from the main grid to operate in either grid-connected or islanded mode supporting customer critical resources as needed. As it has advanced from early prototypes to relatively mature technologies, the microgrid concept has attracted growing attentions from both industry and academia for implementation and deployments [2].

Interfacing between the utility grid and customers, microgrids are effective solutions to provide multiple services to both parties with reduced operating costs [3], [4], [5], improved reliability [6], [7], [8] and enhanced resiliency [9], [10], [11]. Optimal DER sizing to determine energy storage [12], [13], PV [14], [15] and demand response [16], [17], [18] capacities are essential to gain all the benefits microgrids can offer. Chen *et al.* [3] sized energy storage from a cost-benefit analysis with microgrid unit commitment scheduling. A cooperative planning approach was proposed for networked microgrids considering wind and solar [5]. Bahramirad *et al.* [6] formulated a mixed-integer linear program (MILP) to invest energy storage to minimize the operating costs as well as the reliability impacts from generator outages and renewable intermittency. Similarly, community microgrids reliability was taken into account as in [7], [8]. Resiliency against natural disasters is a pressing issue, and Khodaei [9] proposed a decomposition framework to iteratively solve for optimal schedule. A robust approach to consider the renewable, load and market price uncertainty was introduced by the same author in [19].

Compared to the conventional AC microgrids mentioned before, DC microgrids [20], [21] recently are gaining popularity for their direct links with DC generation and loads, such as PV, ES, LED lightning and DC data centers. However, it is still costly to build or upgrade to pure DC microgrids. Hence, hybrid AC/DC microgrids [22] are good alternatives to leverage all the benefits from both AC and DC microgrids without huge construction upgrades. This is particularly appealing to achieve widespread commercial viability. Lotfi *et al.* has compared AC or DC configurations [24] and proposed a hybrid microgrid formulation [23].

Currently, commercial building owners pay expensive electricity bills for both energy usage charges and peak demand charges and follow stringent reliability requirements for the critical loads inside. Converting these buildings into hybrid microgrids with properly sized DER augments the overall economics and reliability. However, existing microgrid sizing approaches either neglect two-part charges with undersized DER or are conservative leading to oversized ones. Full formulation for the hybrid microgrid is missing as well.

This paper proposes a two-stage optimal DER sizing strategy for commercial building hybrid microgrids. It jointly minimizes the investments of energy storage, PV and power electronic devices as well as the expected operating costs which include grid-connected mode energy charges, demand charges, storage degradation costs and islanded mode load shedding penalties. This formulation is then linearized into a computational efficient MILP. With optimized DER capacities, this model achieves economic grid-connected mode and reliable islanded mode performances. Numerical results based on real commercial building load profiles also demonstrate the importance and necessity of simultaneously considering grid-connected mode energy charges and demand charges together with islanded load shedding reliability costs.

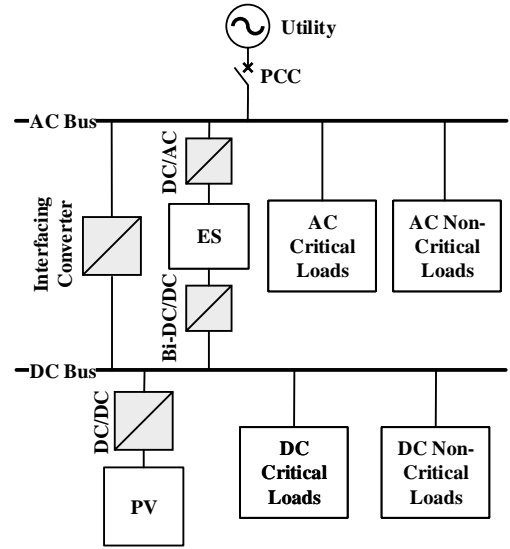

Fig. 1. Hybrid Microgrid Configuration

The rest of the paper is organized as follows. Section II presents the commercial building hybrid AC/DC microgrid configuration and the mathematical formulation of the proposed sizing approach. Section III presents a detailed case study based on the real data collected from a Californian commercial building. The influence of the energy charges, demand charges and critical loads are discussed. Section IV concludes the paper.

## II. MATHEMATICAL FORMULATION

This section provides the mathematical formulation of a DER sizing strategy for a commercial building hybrid AC/DC microgrid. It optimizes the ES, PV and power electronics devices capacities to minimize the combination of investments, grid-connected electricity bills, islanded load shedding penalties and storage degradation.

### A. Commercial Building Hybrid Microgrid Configuration

Fig. 1 illustrates the adopted hybrid microgrid configuration in this paper. The microgrid is designed for a commercial building with critical load resources as data centers. This microgrid includes an AC bus and a DC bus where critical loads are attached. ES is shared by both buses via a DC/AC inverter and a DC/DC converter. The other DER component, PV, is linked to the DC bus by a DC/DC converter for maximum power point tracking (MPPT) and direct power control. Utility grid is connected at the point of common coupling (PCC) for AC bus voltage and frequency control and collecting energy and demand charges from the microgrid (building) owners. Reverse flow is not permitted for utility regulation purpose.

### B. Optimal Sizing Strategy

$$\min \quad obj = c^{inv} + \sum_s \pi_s (c_s^e + c_s^d + c_s^{lcl} + c_s^{lnl} + c_s^{deg}) \quad (1)$$

$$c^{inv} = C^{PV} x^{PV} + C^{ES} x^{ES} + C^{IC} x^{IC} \\ + C^{INV} x^{INV} + C^{CON} x^{CON} \quad (2)$$

$$c_s^e = \sum_t \lambda_t^e p_{s,t}^{grid} \quad (3)$$

$$c_s^d = \lambda^d p_s^{peak} \quad (4)$$

$$c_s^{deg} = \sum_t c^{deg} (dch_{s,t}^{AC} + dch_{s,t}^{DC} + ch_{s,t}^{AC} + ch_{s,t}^{DC}) \quad (5)$$

$$c_s^{lcl} = \sum_t VOLL^{CL} (lcl_{s,t}^{AC} + lcl_{s,t}^{DC}) \quad (6)$$

$$c_s^{lnl} = \sum_t VOLL^{NL} (lnl_{s,t}^{AC} + lnl_{s,t}^{DC}) \quad (7)$$

$$dch_{s,t}^{AC} \eta^{INV} - ch_{s,t}^{AC}/\eta^{INV} + p_{s,t}^{grid} = \\ f_{s,t}^{AC} + CL_{s,t}^{AC} + NL_{s,t}^{AC} \quad (8)$$

$$(dch_{s,t}^{DC} + v_{s,t}^{DC})\eta^{CON} - ch_{s,t}^{DC}/\eta^{CON} = \\ f_{s,t}^{DCout} - f_{s,t}^{DCin} + CL_{s,t}^{DC} + NL_{s,t}^{DC} \quad (9)$$

$$f_{s,t}^{AC} = f_{s,t}^{DCin}/\eta^{IC} - f_{s,t}^{DCout} \eta^{IC} \quad (10)$$

$$0 \leq f_{s,t}^{DCin} \leq z_{s,t}^{DC} M \quad (11)$$

$$0 \leq f_{s,t}^{DCout} \leq (1 - z_{s,t}^{DC}) M \quad (12)$$

$$dch_{s,t}^{AC} + dch_{s,t}^{DC} \leq x^{ES} y_{s,t}^{ES} \quad (13)$$

$$ch_{s,t}^{AC} + ch_{s,t}^{DC} \leq x^{ES} (1 - y_{s,t}^{ES}) \quad (14)$$

$$\alpha^{min} \rho^{EP} x_{ES} \leq soc_{s,t} \leq \alpha^{max} \rho^{EP} x_{ES} \quad (15)$$

$$soc_{s,t} = soc_{s,t-1} + (ch_{s,t}^{AC} + ch_{s,t}^{DC})\eta^{ch} \\ - (dch_{s,t}^{AC} + dch_{s,t}^{DC})/\eta^{dch} \quad (16)$$

$$v_{s,t}^{DC} \leq V_{s,t} x^{PV} \quad (17)$$

$$0 \leq p_{s,t}^{grid} \leq p_s^{peak} \quad (18)$$

$$p_s^{peak} \leq P^{Tariff} \quad (19)$$

$$dchi_{s,t}^{AC} \eta^{INV} = fi_{s,t}^{AC} + CL_{s,t}^{AC} - \\ lcl_{s,t}^{AC} + NL_{s,t}^{AC} - lnl_{s,t}^{AC} \quad (20)$$

$$(vi_{s,t}^{DC} + dchi_{s,t}^{DC})\eta^{CON} + fi_{s,t}^{DCin} - fi_{s,t}^{DCout} = \\ CL_{s,t}^{DC} - lcl_{s,t}^{DC} + NL_{s,t}^{DC} - lnl_{s,t}^{DC} \quad (21)$$

$$0 \leq lcl_{s,t}^{AC} \leq CL_{s,t}^{AC} \quad (22)$$

$$0 \leq lnl_{s,t}^{AC} \leq NL_{s,t}^{AC} \quad (23)$$

$$0 \leq lcl_{s,t}^{DC} \leq CL_{s,t}^{DC} \quad (24)$$

$$0 \leq lnl_{s,t}^{DC} \leq NL_{s,t}^{DC} \quad (25)$$

$$fi_{s,t}^{AC} = fi_{s,t}^{DCin}/\eta^{IC} - fi_{s,t}^{DCout} \eta^{IC} \quad (26)$$

$$0 \leq fi_{s,t}^{DCin} \leq zi_{s,t}^{DC} M \quad (27)$$

$$0 \leq fi_{s,t}^{DCout} \leq (1 - zi_{s,t}^{DC}) M \quad (28)$$

$$vi_{s,t}^{DC} \leq V_{s,t} x^{PV} \quad (29)$$

$$dchi_{s,t}^{AC} + dchi_{s,t}^{DC} \leq x^{ES} \quad (30)$$

$$dchi_{s,t}^{AC} + dchi_{s,t}^{DC} \leq soc_{s,t-1} \quad (31)$$

The objective function (1) minimizes the annualized microgrid investments and expected operating and reliability costs. Equation (2) describes the purchase costs for PV, ES, bi-directional interfacing converters (IC), AC/DC inverters and DC/DC converters. Utility bills for energy and demand charges are calculated in (3) and (4). Grid-connected mode storage degradation is also accounted in (5) to avoid frequent deep cycle charging and discharging. Critical and non-critical load shedding during the islanded mode are penalized in (6)–(7).

Grid-connected mode operation are modeled in (8)–(19). Equations (8) and (9) enforce the power balance for the AC and DC buses where load shedding is not allowed. The power flow through the interfacing converter is formulated in constraints (10)–(12). Due to the bidirectional conversion loss, binary variables $z_{s,t}^{DC}$ are used to indicate the flow direction with big-M constraints as in [12], [25]. ES charging, discharging, state-of-charge (SoC) limits and SoC transitions are described in (13)–(16). To avoid simultaneous charging/discharging, ES discharging states are also optimized with binary variables $y_{s,t}^{ES}$. With MPPT and direct power control, PV generation should not exceed the maximum available output as in (17). Following tariff plan, peak demand for the net load is calculated and limited in (18) and (19).

Islanded mode is one of the most salient features of microgrids, which enables a disconnected service from the main feeder. During the islanded mode, the proposed hybrid microgrid needs to minimize the involuntary load shedding especially for the critical loads. One-hour-long backup is required here for every selected operating period. Constraints (20)–(31) represent this reliable and resilient operation where energy are reserved in ES to provide necessary support along with PV.

With site physical space limits, constraints (32) and (33) indicate the maximum allowable PV and ES sizes. Power electronic capacities are then determined based on the grid-connected and islanded operation as in (34)–(41).

$$x^{PV} \leq PV^{max} \quad (32)$$

$$x^{ES} \leq ES^{max} \quad (33)$$

$$x^{INV} \geq dch_{s,t}^{AC} + ch_{s,t}^{AC}/\eta^{INV} \quad (34)$$

$$x^{INV} \geq dchi_{s,t}^{AC} \quad (35)$$

$$x^{CON} \geq x^{PV} + dch_{s,t}^{DC} + ch_{s,t}^{DC}/\eta^{CON} \quad (36)$$

$$x^{CON} \geq x^{PV} + dchi_{s,t}^{DC} \quad (37)$$

$$x^{IC} \geq f_{s,t}^{DCin}/\eta^{IC} \quad (38)$$

$$x^{IC} \geq f_{s,t}^{DCout} \quad (39)$$

$$x^{IC} \geq fi_{s,t}^{DCin}/\eta^{IC} \quad (40)$$

$$x^{IC} \geq fi_{s,t}^{DCout} \quad (41)$$

To deal with the nonlinear product $x^{ES} y_{s,t}^{ES}$, auxiliary variables are introduced in (42) to replace nonlinear constraints in (13) and (14) with the exact linear reformulations (43)–(47).

$$x^{ES} y_{s,t}^{ES} = u_{s,t}^{ES} \quad (42)$$

$$u_{s,t}^{ES} = x^{ES} - \kappa_{s,t}^{ES} \quad (43)$$
$$0 \leq u_{s,t}^{ES} \leq My_{s,t}^{ES} \quad (44)$$
$$0 \leq \kappa_{s,t}^{ES} \leq M(1 - y_{s,t}^{ES}) \quad (45)$$
$$dch_{s,t}^{AC} + dch_{s,t}^{DC} \leq u_{s,t}^{ES} \quad (46)$$
$$ch_{s,t}^{AC} + ch_{s,t}^{DC} \leq x^{ES} - u_{s,t}^{ES} \quad (47)$$

## III. CASE STUDY

### A. Simulation Setup

The proposed hybrid microgrid sizing approach has been tested with the load data from a commercial building located at California, USA. Critical loads include the AC/DC loads in data centers, AC cooling and heating as well as necessary AC/DC lighting and appliances. The peak load at this building is at 846 kW. The PV profiles are based on NREL PVWATTTS Dataset [26]. For planning purpose, the complete annual load could be used to capture the temporal evolution. However, since the load data are similar with a limited number of patterns, scenario reduction technique [27] is applied here to select typical days to preserve accurate load variations as well as to improve computational efficiency. In this simulation, 6 days are chosen to represent different load shapes in different months, and PV profiles are selected accordingly.

The planning horizon for this microgrid is considered to be 10 years with 10% annual discount rates. Annualized investment costs and the technical parameters are listed in the Table I. A two-hour Lithium-ion battery is adopted here. This microgrid participates the PG&E E19 tariff plan where the peak load should be below 1000 kW. The load shedding penalties for critical and non-critical loads are selected as $ 3,000 /kWh and $ 500 /kWh, respectively.

TABLE I. INVESTMENT PARAMETERS

| | Annualized Investments | Allowable Capacities | Technical Parameters |
|---|---|---|---|
| PV | $ 108 /kW | 400 kW | N/A |
| ES | $ 424 /kW | 350 kW | $\eta^{dch} = \eta^{ch} = 0.93$ |
| DC/AC Inverter | $ 6.5 /kW | N/A | $\eta^{INV} = 0.96$ |
| DC/DC Converter | $ 4.3 /kW | N/A | $\eta^{CON} = 0.98$ |
| Interfacing Converter | $ 8.1 /kW | N/A | $\eta^{IC} = 0.96$ |

All simulations were carried out in Julia using the JuMP package [28] and CPLEX solver on an Intel Core 4.00 GHz processor with 12GB of RAM. The MIP optimality gap was set to 0.01%. The approximate running time for each simulation is less than 1s.

### B. Numerical Results

Four cases are studied in this paper. 1) *Case 0*: Base case with no DER; 2) *Case 1*: PV only deployment; 3) *Case 2*: ES only deployment; 4) *Case 3*: full DER deployment. The planning results are summarized in the Table II. From the table, PV and ES installations always hit the physical sizing limits to gain economic and reliability benefits for investment recovery. However, the benefits from these two DER are not quite the same. PV mainly leads to electricity bill reduction and ES contributes to the reliability enhancement.

PV shapes are highly correlated with the building loads, so the PV outputs can directly shave the peak to reduce the energy charges and demand charges at the same time. Compared with the base case when no DER is installed, *Case 1* PV only case achieves 27.1% energy charges savings, 17.5% demand charges savings and 25.6% total bill payment savings. On the contrary, *Case 2* only obtains 1.2% energy charges savings, 11.6% demand charges savings and 2.7% total bill payment savings, which is surprisingly low. The reason is that unlike PV, ES cannot generate power from its own, and it has to charge first from either PV or the grid. The price differences between the peak hours and off-peak hours plus demand charge reductions have to compensate for the storage degradation costs, the charging/discharging loss and the energy conversion loss through inverters, converts or IC. This requirement substantially hinges ES peak shaving behaviors especially when the PV is absent. Combining the advantages from both PV and ES, *Case 3* DER case demonstrates 28.5% energy charges savings, 35.3% demand charges savings and 29.5 % total savings. Interestingly, these cost savings numbers are all larger than the sum of the cost savings from each single deployment (*e.g.* 28.5% > 27.1% + 1.2%) even though the exceeding amount is not large. This confirms the benefits to the microgrids from the DER joint operation.

TABLE II. INVESTMENT RESULTS

| | *Case 0*: Base Case | *Case 1*: PV only | *Case 2*: ES only | *Case 3*: DER |
|---|---|---|---|---|
| PV (kW) | 0 | 400 | 0 | 400 |
| ES (kW) | 0 | 0 | 350 | 350 |
| DC/AC inverter (kW) | 0 | 0 | 211 | 350 |
| DC/DC converter (kW) | 0 | 400 | 276 | 624 |
| Interfacing converter (kW) | 353 | 228 | 295 | 199 |
| Energy charges ($ $10^4$) | 74.59 | 54.34 | 73.68 | 53.34 |
| Demand charges ($ $10^4$) | 12.70 | 10.55 | 11.23 | 8.22 |
| Total payment ($ $10^4$) | 87.29 | 64.89 | 84.91 | 61.56 |

Fig. 2 presents the load curtailments for these four cases. When there is no DER in the system, inconvenient load shedding occurs for both critical loads and non-critical loads during the islanded mode. *Case 1* PV only case reduces the curtailment to a extent during the daytime, but it cannot contribute at evening or night when there is no sunlight. Even though *Case 2* ES only case does not help much on the bill reduction, the reliability improvement is significant. The critical loads are fully covered and non-critical loads are well preserved as well. Again, *Case 3* DER case leads to the most reliable case which improves system emergency support.

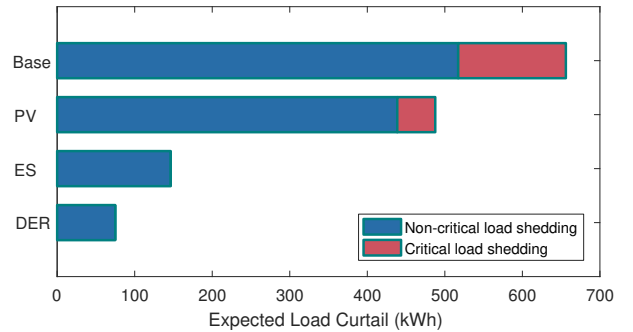

Fig. 2. Curtailment for critical load and non-critical load

Fig. 3 shows the net load profiles before/after DER deployment and the corresponding DER outputs. The upper plot presents a peak day scenario. Before installing PV and ES, the peak demand is 846 kW, whereas it is only 454 kW after DER integration. This reduced peak is mainly contributed from PV generation. ES also discharge to help reducing the peak at noon time. As PV only generates during the daytime, the energy reduction at rest of the day is not significant. The lower plot displays a normal day scenario. At this particular day, PV output can fully covers the commercial building load and leads to a "Duck Curve" shaped load. In addition, due to the islanded mode constraints, ES has scheduled to store adequate energy to hedge against potential reliability and resiliency events.

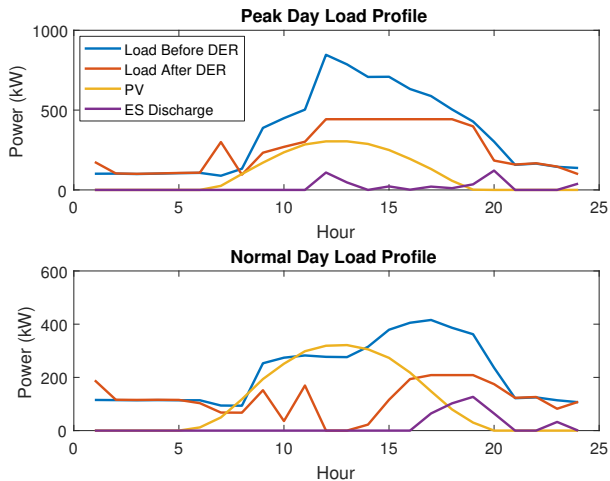

Fig. 3. Load profiles before and after DER integration

## IV. Conclusion

In this paper, a new MILP-based commercial building hybrid AC/DC microgrid sizing model is proposed to leverage the benefits from both AC and DC microgrids. Energy storage, photovoltaics and power electronic devices are optimally sized to coordinate between grid-connected mode economics and islanded mode reliability. The model integrates the electricity bill reductions into the optimization scheme and enables a minimized involuntary load shedding to secure critical loads. Numerical results from a commercial building load demonstrate the effectiveness of the proposed model on economics and reliability enhancement with DER in the hybrid microgrid.